\documentclass[12pt]{article}
\usepackage{latexsym}
\usepackage[latin1]{inputenc}
\usepackage{epsfig}
\usepackage{array}
\usepackage{amssymb}
\usepackage{amsmath}
\usepackage{amsthm}
\newcommand{\remark}[1]{}
 
\def\P{\operatorname{\bf P}}
\newcommand{\bit}{\begin{itemize}}
\newcommand{\eit}{\end{itemize}}
\newcommand{\bc}{\begin{center}}
\newcommand{\ec}{\end{center}}
\newcommand\G{\mathbf{G}}
\newcommand\GW{\mathbf{G}}
\newcommand\T{\mathbf{T}}
\newcommand\N{\mathbf{N}}

\newcommand\D{\mathbf{D}}

\newcommand\var{\mathbf{\sigma}}
\newcommand\tiempo{{\sl time}}
\begin{document}
\title{Counting Knight's Tours through the \\
Randomized Warnsdorff Rule}
\author{H\'ector Cancela\footnote{Instituto de Computaci\'on, Facultad de Ingenier\'{\i}a, Universidad de la República. Address:
Julio Herrera y Reissig 565, Montevideo, Uruguay.
Código Postal: 11300.
Casilla de Correo N. 30.
e-mail: cancela@fing.edu.uy}
\quad 
and
\quad 
Ernesto Mordecki\footnote{Centro de Matem\'atica, Facultad de Ciencias, 
Universidad de la República. Address: Facultad de Ciencias,
Igu\'a 4225, Código Postal: 11400, Montevideo, Uruguay.
e-mail: mordecki@cmat.edu.uy}}
\maketitle
\begin{abstract}
We give an estimate of the number of geometrically distinct open
tours $\G$ for a knight on a chessboard. We use a randomization of 
Warnsdorff rule to implement importance sampling in a backtracking
scheme, correcting the observed bias of the original rule, according to the proposed 
principle that
``most solutions follow Warnsdorff rule most of the time''.
 After some experiments in order to test this principle, 
and to calibrate a parameter, interpreted
as a distance of a general solution from a Warnsdorff solution,
we conjecture that $\G=1.22\times 10^{15}$.   
\end{abstract}
\newpage
\section{Introduction}
In this paper we give an estimate of the number of open Knight Tours in a chessboard, 
based on the heuristics given by the classical Warnsdorff rule, 
used to implement importance sampling 
in a backtracking algorithm, following [Knuth 1975].
The historical development of the Knight Tour 
problem is beyond the possibilities of the present work,
and the interested reader is suggested to visit George Jelliss's web page
and the references therein [Jelliss 2003], from where we borrow our 
notation.
Let us define precisely the quantities we are interested in:
\begin{itemize}
\item
$\G$ is the number of $\G$eometrically distinct open tours.
\item
$\T$ is the number of open $\T$our diagrams, by rotation and symmetry 
$\T = 8\,\G$.
\item
$\N$ is the number of open tour $\N$umberings: 
$\N = 2\,\T = 16\,\G$ since each can be numbered from either end.
\item
$\D$ is the number of closed tour $\D$iagrams.
\end{itemize} 
These parameters refer to the standard $8\times 8$ chessboard.
When other sizes of (squared) chessboards are considered, we include a subscript in our notation,
for instance $\N_6$ is the number of open tour numberings in a $6\times 6$ chessboard.

In the terms just introduced our purpose is to give an estimate of $\G$.
A relevant exact computation related to this problem
was carried out independently in [McKay 1997] and in [Wegener 2000], 
who computed the number
$\D=13,267,364,410,532$.

\section{Importance sampling through Warnsdorff rule}
\label{iswr}

Warnsdorff rule was  proposed in 1823 in order to find open tours in a chessboard.
It states: ``at each step,
choose the square with the lowest number of possible continuations;
if a tie occurs it may be broken arbitrarily'' [Warnsdorff 1823]. 
This heuristic rule is incredibly effective to find solutions.
Only 150 years after being formulated, with the help of a computer, 
it was reported that the rule can fail, in respect to its second part [Guik 1983].
Modern numerical experiments showed that this rule is \emph{biased},
in the sense that it founds solutions in an 
extremely small subset of the set of solutions, fact that we confirm.

Our departure point in order to estimate the total amount of open tours, 
is the principle that ``most solutions follow Warnsdorff rule most of the time'', 
i.e. a typical solution is a concatenation of short paths satisfying 
the rule within them, that, when concatenated violate the rule.
According to this principle we decompose the number $\G$ as the sum
$\G=\GW_0+\cdots+\GW_{62}$
where we denote by $\GW_k$ ($k=0,1,\dots,62$)
the number of geometrically distinct knight tours that violates the rule 
exactly $k$ times. In particular $\GW_0$ is the number of knight tours
that follows Warnsdorff rule 
(including any of the possibilities in case of tie). 
In section \ref{w} we provide some empirical evidence supporting this principle,
giving an estimation of $\GW_k$ for different values of $k$.
We obtain the largest number when $k=13$,
and significant numbers in a range from $k=5$ to $k={20}$.

In consequence, following  [Knuth 1975] and these considerations, 
we implement a randomized backtracking algorithm, that we call the
``Randomized Warnsdorff Rule''.
\begin{figure}[ht]
\begin{center}\label{tablerito}
\epsfig{file=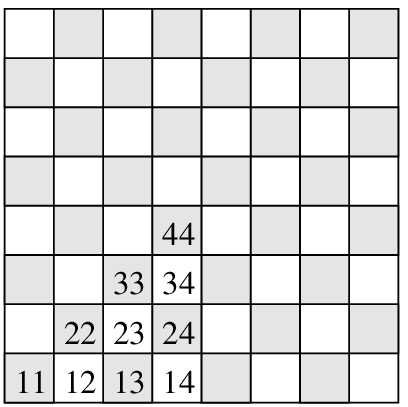}
\caption{Relevant initial squares}
\end{center}
\end{figure}
We begin by fixing the ten geometrically different initial  squares (see Figure \ref{tablerito})
and, at step $n$,
we assign probabilities $p(n,\alpha,i(n),j)$ to a jump from the current square $i(n)$ to each connected square $j$,
proportional to $N_{j}^{-\alpha}$ where $N_{j}$ 
is the number of free accessible squares from $j$
and $\alpha$ is an \emph{Importance Sampling} (IS) parameter.
The estimation of the number of tours for each run is
$$
\begin{cases}
\prod_{n=1}^{63}\Big(p\big(n,\alpha,i(n),i(n+1)\big)^{-1}\Big) & \text{ when the algorithm finds the solution $\{i(n)\}$,}\\
0&\text{ when the algorithm does not complete a tour.}
\end{cases}
$$
This is an \emph{unbiased} estimate of the number of tours beginning at $i(1)$, 
for all  positive probability assignations [Knuth 1975].
The problem is then to choose the probabilities 
$p\big(n,\alpha,i,j\big)$
in order to achieve
the best possible \emph{variance reduction}, and here we 
use Warnsdorff heuristics.

Two important cases of our scheme are $\alpha=0$, corresponding to plain
Monte Carlo simulation in the graph, and $\alpha=\infty$,
that gives maximum (equal) probability to the nodes with minor amount
of free nodes, and null probability to the others, i.e. the original
Warnsdorff rule with equiprobability in case of tie, 
giving a biased estimate. 

The first step of our work consists in calibrating the parameter $\alpha$,
that measures the distance form a generic solution from a Warnsdorff solution, 
running our randomized Warnsdorff algorithm
to estimate three known quantities, $\N_5$, $\D_6$ and $\D$. 
It is important to note (and it has been noticed in the literature)
that the optimal $\alpha$ is sensible to the size of the chessboard.
We report the results in Section \ref{tuning}

The second step, once we calibrate $\alpha=1.5$, 
consists in running a relatively long sampling experiment to obtain our final estimate of $\G$,
presented in Section \ref{results}. 
We use this run also to compute the results presented in Section \ref{w}.

\section{Parameter tuning}
\label{tuning}

In this section we present the results of some numerical experiments that were 
designed in order to validate the proposed method and at the same time to 
tune $\alpha$, the importance sampling parameter.

We estimate three different known quantities using the randomized Warnsdorff rule:
$\N_5=1728$, obtained in [Planck 1908] and [Planck 1909],
$\D_6=98626$ obtained in [Duby 1964],
and the McKay-Wegener number $\D$.

To define the experimental design, it is necessary to fix a certain number of 
parameters, in particular the sampling size for the Monte Carlo experiments, 
the number of replications with different seeds, and the importance sampling 
parameter $\alpha$. The two first parameters 
were fixed taking into account some preliminary experiments and the need to 
complete work in a reasonable computing time; the sampling size was then 
fixed to $5\times 10^5$ for the $5\times 5$ board and to $10^6$ for the 
other two cases, and 21 replications were made with different seeds. Even 
if $\alpha$ is a continuous parameter, it is
necessary to decide on a (small) set of values for employing in the experiments. 
It is important to include the case where there is no IS at all ($\alpha=0$), 
as well as a relatively wide range of values
where the method has different behaviors. We also included $\alpha=-1$, 
which corresponds to an IS scheme which favors tours 
``far away from Warnsdorff rule''. 
Based on the same preliminary experiments, we choose to take 
$\alpha$ in the set 
$\{-1, 0, 0.5, 1, 1.5, 2, 2.5, 3, 4, 5, 8, 10\}$.
The total number of runs is then 756 (three problems, 12 different values for $\alpha$, 
21 independent replications with different seeds for each of these). 
The Monte Carlo method was programmed using C 
(compiled with gcc 3.3.3 and 3.4.3), and all the 
experiments were run on an personal computer (CPU 1.70GHz, RAM 224508 kB) running Linux.
Each experiment took in the order of 100 CPU seconds on the 
$5\times 5$ board, 2000 CPU seconds on the $6\times 6$ 
board and 15000 seconds on $8\times 8$ boards (as reported
by the {\tt times()} function of the {\tt ``sys/times.h''} library).
\begin{table}[h]
\centering
\small
\begin{tabular}{|l|c|>{\bf}c|c|c|c|}
\hline
$\alpha$   & $\min_i\{\N_5(i)\}$ & $\widehat{\N_5}$ & $\max_i\{\N_5(i)\}$ & $\var$ & \tiempo \\
\hline 
-1    &   314  &   2120      &   8561    &    2168        &       75.3   \\
 0    &  1578  &   1732      &   1887    &      83        &       56.0   \\
0.5   &  1663  &   1716      &   1789    &      28        &       129.6  \\
1     &  1693  &   1729      &   1758    &      13        &       100.5  \\
1.5   &  1714  &   1729      &   1742    &       8.0      &       138.2  \\
2     &  1720  &   1729      &   1741    &       6.2      &       118.5  \\
\hline
2.5   &  1722  &   1728      &  1738     &       6.0      &       153.8  \\
\hline
3     &  1712  &   1728      &   1740    &       6.8      &       118.7  \\
4     &  1702  &   1726      &   1752    &      12        &       128.7  \\
5     &  1686  &   1732      &   1795    &      27        &       129.0  \\
8     &  1279  &   1676      &   2589    &     339        &       133.8  \\
10    &  1002  &   1866      &   11310   &    1871        &       124.4  \\
\hline        
\end{tabular}
\caption{Quality of estimations of $\N_5$={\bf 1728} 
as function of the IS parameter $\alpha$.}
\label{tab5x5}
\end{table}
A summary of the results obtained is shown in Tables \ref{tab5x5}, \ref{tab6x6}, and 
\ref{tab8x8} respectively for the three problems studied. Given a problem, for instance,
the estimation of $\N_5$, the number of open tours numberings in a $5\times 5$
chessboard, and a value for $\alpha$, 
we have $i=1,\dots,21$ independent replications, 
each returning an estimation numbered through
$\N_5(i)$, and an estimation $\var(i)$ of the 
standard deviation (i.e. the square root of the variance)
of this $i$-th estimation. 
The final estimation, the mean in $i$, is denoted $\widehat{\N}_5$. 
Each table has five columns, 
the first one corresponding to the $\alpha$ values, and the 
following ones showing the minimum, the estimation, and maximum values for the 21 estimations,  
an estimation $\var$ of the standard deviation of our estimate, 
and the mean of the execution times. 
The same analysis is performed with $\D_6$ and $\D$.

It comes as no surprise to observe that the first two problems are rather easy, and a 
wide range of $\alpha$ values result in estimations quite near to the exact value. 
It is clear that the variance of the estimation depends strongly on $\alpha$; in both 
cases the best value is observed when $\alpha=2$ 
(but  $\alpha=1.5$ and  $\alpha=2.5$ come very close). High values for $\alpha$ seem to 
greatly affect the precision of the method (this can be specially be seen in the last 
two lines for the $6\times 6$ board, where the estimations are quite bad). 
Also, $\alpha=-1$ leads to inconsistent results.

The standard (no importance sampling) Monte Carlo results can be seen in the $\alpha=0$ 
line. In both cases the standard method gives quite 
consistent estimates for the number of tours, but using Importance Sampling we can 
attain a variance orders of magnitude smaller. 

From these experimental results, it can be seen that in the three
cases under study the proposed method works well and is more 
robust and more efficient than the standard Monte Carlo, which in
some cases does not give meaningful estimations of the desired parameters. 
Fixing $\alpha$ can be a problem, as the method is
clearly dependent on this parameter, but fortunately it seems as if the 
choice of $\alpha$ around $1.5$ or $2$ is rather robust, as those
are the values with best results in these three problems.
As we want to estimate in the bigger chessboard, we decide ourselves for $\alpha=1.5$.
\begin{table}[h]
\centering
\small
\begin{tabular}{|l|c|>{\bf}c|c|c|c|}
\hline
$\alpha$   & $\min_i\{\D_6(i)\}$ & $\widehat{\D_6}$ & $\max_i\{\D_6(i)\}$ & $\var$ & \tiempo \\
\hline
-1& 0   &       0       &       0    &          0       &       1.142E+03       \\ 
0&
7444    &       10113   &       13920    &         1387     &       9.207E+02       \\
0.5&                                                                    
9324    &       9945    &       10547    &          220     &       2.076E+03       \\
1&
9718    &       9861    &       9960    &         64       &       1.649E+03       \\
1.5&
9816    &       9861    &       9925    &         31       &       2.312E+03       \\
\hline
2&
9808    &       9862    &       9887    &         24       &       2.002E+03       \\
\hline
2.5&
9814    &       9871    &       9906    &         26       &       2.783E+03       \\
3&
9799    &       9866    &       9929    &         37       &       2.162E+03       \\
4&
9690    &       9896    &       10521    &         191      &       2.416E+03       \\
5&
9231    &       10074   &       14568   &        1105      &       2.444E+03       \\
8&
5252    &       8234    &       17429    &        2610      &       2.761E+03       \\
10&
3162    &       5866    &       25270    &        4888      &       2.575E+03       \\ 
\hline
\end{tabular}
\caption{Quality of estimations of $\D_6= {\bf 9862}$ as function of the IS parameter $\alpha$.}
\label{tab6x6}
\end{table}
Computing the number of closed tour diagrams on the $8\times 8$ board is a more demanding 
problem, and the standard (no importance sampling) method can not find good estimates (nor can  do it the method with $\alpha=-1$). 
In general, the IS estimates quality depends heavily on $\alpha$, and only when 
$0.5\leq \alpha\leq 3$ we find estimation values that approach the exact value. 
Here  the setting $\alpha=1.5$ leads to the lowest variance.
\begin{table}[h]
\centering
\small
\begin{tabular}{|l|c|>{\bf}c|c|c|c|}
\hline
$\alpha$   & $\min_i\{\D(i)\}$ & $\widehat{\D}$ & $\max_i\{\D(i)\}$ & $\var$ & \tiempo \\
\hline 
-1 &0           &       0               &               0       &       0               &       8.517E+03       \\
0&                                      
0               &       2.992E+13       &               3.999E+13       &      8.82E+13        &       5.912E+03       \\
0.5&
1.061E+13       &       1.386E+13       &       1.692E+13       &      1.56E+12       &       1.339E+04       \\
1&
1.253E+13       &       1.320E+13       &       1.391E+13       &      3.48E+11       &       1.038E+04       \\
\hline
1.5&
1.280E+13       &       1.318E+13       &       1.398E+13       &      2.36E+11        &       1.507E+04       \\
\hline
2&
1.248E+13       &       1.324E+13       &       1.463E+13       &      4.17E+11       &       1.471E+04       \\
2.5&
1.225E+13       &       1.315E+13       &       1.534E+13       &      6.65E+11       &       1.839E+04       \\
3&
1.104E+13       &       1.327E+13       &       2.456E+13       &      3.12E+12       &       1.554E+04       \\
4&
7.107E+12       &       1.122E+13       &       2.413E+13       &      4.35E+12       &       1.717E+04       \\
5&
3.228E+12       &       7.029E+12       &       3.658E+13       &      7.10E+12        &       1.729E+04       \\
8&
7.694E+10       &       3.472E+11       &       1.571E+12       &       4.31E+11      &       1.763E+04       \\
10&
1.013E+10       &       1.687E+11       &       2.897E+12       &       6.28E+11       &       2.165E+04       \\
\hline                                                                \end{tabular}
\caption{Quality of estimations of $\D= {\bf 13,267,364,410,532}$ as a function of $\alpha$.}
\label{tab8x8}
\end{table}
\section{On the set of open tours}
\label{w}
\begin{figure}[ht]
\begin{center}
\epsfig{file=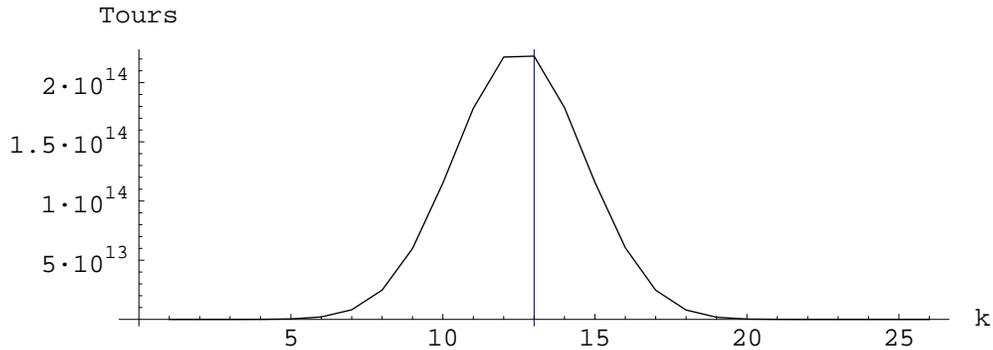}
\caption{Estimated number of tours $\GW_k$  for $k=0,\dots,25$ ($N=10^9$).}
\label{histogram}
\end{center}
\end{figure}

In this section we present some empirical evidence supporting our departure
point, that is, that ``most solutions follow Warnsdorff rule most of the time''.
Firstly, we present in Figure \ref{histogram} 
a histogram where we show the estimations of $\GW_k$ obtained in the final 
simulation (described in next section), in the case
and $\alpha=1.5$, with $10^9$ replications.

In the Table \ref{numbers} we present the corresponding values
plotted in the Figure \ref{histogram}. It is concluded that
values from $\G_5$ to $\G_{20}$ concentrate the 99.9892\% (practically all) solutions,
and also that $\G_0<<\G$, i.e. the number of Warnsdorff solutions is several orders
smaller than the number of solutions. 
\begin{table}[h]
\centering
\begin{tabular}{|c|c|c|c|c|c|c|c|}
\hline
$k$ & $\GW_k$    & $k$    & $\GW_k$       &  $k$      & $\GW_k$   &  $k$       & $\GW_k$     \\
\hline
0   & 1.066E+05  &   6    & 2.084E+12     &   12   &   2.216E+14  &   18       & 7.920E+12   \\ 
1   & 2.357E+07  &   7    & 8.110E+12     &   13   &   2.224E+14  &   19       & 1.845E+12   \\
2   & 5.095E+08  &   8    & 2.475E+13     &   14   &   1.789E+14  &   20       & 3.063E+11   \\
3   & 6.935E+09  &   9    &  5.980E+13    &   15   &   1.158E+14  &   21       & 6.020E+10   \\
4   & 6.300E+10  &  10    &  1.152E+14    &   16   &   6.050E+13  &   22       & 1.121E+09   \\
5   & 4.160E+11  &  11    &  1.783E+14    &   17   &   2.462E+13  &   23       & 0           \\
\hline
\end{tabular}
\caption{Estimations of $\GW_k$ corresponding to Figure \ref{histogram}}
\label{numbers}
\end{table}
\section{On the number $\G$ of open tours}
\label{results}
In this section we present our main results, corresponding to the final numerical experiment.
In this case we perform the randomized Warnsdorff algorithm
(as explained in \ref{iswr}),
with $\alpha=1.5$, as suggested by the
experiments exposed in section \ref{tuning}.
We choose a sample size corresponding to $10^{9}$ size,
as this seems to provide a reasonable result
in an acceptable computational time. In our case
CPU time was $6.5\times 10^{5}$ seconds
$\sim 7.5$ days, on a similar personal computer.

In our runs we distinguish according to
the 10 different initial squares (see Figure \ref{tablerito}), 
and, in the four diagonal squares 
we estimate half of the solutions, based on the symmetry.
The dynamic nature of Warnsdorff rule does not seem to
give a simple way of exploiting the ``go and return'' symmetry, and,
in consequence, our crude algorithm gives an estimation of $2 \G$.
To have an idea of the error of our method
we counted also the closed tours, obtaining the estimation
$\widehat{\D}$=1.326390E+13 (compare with the exact value $\D$=1.3267364410532E+13) with
a standard deviation $\sigma$=3.03E+10.

Our final results are presented in Table \ref{main}.
\begin{table}[h]
\centering
\begin{tabular}{|c|c|c|}
\hline
$\widehat{\G}$   & $\G$                & $\var$        \\
\hline
1.222801E+15 & Unknown             & 8.26E+11  \\  
\hline
\end{tabular}
\label{main}
\caption{Estimation of $\G$ (sample size 10E+9).}
\end{table}
%
%
Based on these final results we construct a confidence interval for $\G$. 
Taking into account that for a normal random variable $X$ with mean $\G$ 
and variance $\sigma$ we have
$
\P(\G-3\sigma<X<\G+3\sigma)=0.997,
$
we conclude that 
$
\G\in [1.220,1.225]\times 10^{15}
$
and from there our main conclusion $\G=1.22\times 10^{15}$.

\end{document}